\documentclass[11pt,reqno]{amsart}
\usepackage{amsmath,amssymb,amsthm,mathrsfs,dsfont}
\usepackage{CJK}
\usepackage{amsmath}
\usepackage{dsfont}
\usepackage{mathrsfs}
\usepackage{amsmath,amssymb}
\usepackage{amsfonts}
\usepackage{hyperref}
\usepackage{amsthm}
\usepackage{graphicx}
\usepackage{subfigure}
\usepackage{xcolor}
\usepackage{overpic}

\newtheorem{theorem}{Theorem}[section]
\newtheorem{lemma}{Lemma}[section]
\newtheorem{definition}{Definition}[section]
\newtheorem{proposition}{Proposition}[section]
\numberwithin{equation}{section}
\newtheorem{remark}{Remark}[section]

\newfont{\bb}{msbm10 at 12pt}

\newcommand{\bal}{\begin{aligned}}      \newcommand{\eal}{\end{aligned}}
\newcommand{\ba}{\begin{array}}      \newcommand{\ea}{\end{array}}
\newcommand{\bc}{\begin{center}}     \newcommand{\ec}{\end{center}}
\newcommand{\be}{\begin{enumerate}}  \newcommand{\ee}{\end{enumerate}}
\newcommand{\beq}{\begin{eqnarray}}  \newcommand{\eeq}{\end{eqnarray}}
\newcommand{\beQ}{\begin{eqnarray*}} \newcommand{\eeQ}{\end{eqnarray*}}
\newcommand{\bi}{\begin{itemize}}    \newcommand{\ei}{\end{itemize}}
\newcommand{\bt}{\begin{tabular}}    \newcommand{\et}{\end{tabular}}
\newcommand{\bdm}{\begin{displaymath}} \newcommand{\edm}{\end{displaymath}}

\def\qed{\hfill{Q.E.D.}\smallskip}

\begin{document}

\title{\bf A new proof of Bowers-Stephenson conjecture}
\author{Xu Xu}
\maketitle

\begin{abstract}
Inversive distance circle packing on surfaces was introduced by Bowers-Stephenson \cite{BS} as
a generalization of Thurston's circle packing and conjectured to be rigid.
The infinitesimal and global rigidity of circle packing with nonnegative inversive distance were
proved by Guo \cite{Guo} and Luo \cite{L} respectively.
The author  \cite{X2} proved the global rigidity of circle packing with inversive distance in $(-1,+\infty)$.
In this paper, we give a new variational proof of the Bowers-Stephenson conjecture
for inversive distance in $(-1,+\infty)$,
which simplifies the proofs in \cite{Guo,L,X2} and could be generalized to three dimensional case.
\end{abstract}

\section{Introduction}

In the study of hyperbolic structure on 3-dimensional manifolds,
Thurston \cite{T1} introduced circle packing with non-obtuse intersection angles on surfaces,
which generalized the circle packing studied
by Andreev \cite{An1,An2} and Koebe \cite{K1}.
Thurston proved the Andreev-Thurston theorem, which includes the existence part and the rigidity part.
The Andreev-Thurston rigidity theorem states that the circle packing is globally determined by the discrete curvature on the triangulated surface,
which is defined to be $2\pi$ less the cone angle at a vertex.
Recently, Andreev-Thurston theorem was generalized by Zhou \cite{Z} to the case of obtuse angles.
For a proof of Andreev-Thurston Theorem, see \cite{CL, DV, H, MR, S, T1,Z}.

Inversive distance circle packing was introduced by Bowers-Stephenson \cite{BS} as a generalization of
Thurston's circle packing on surfaces, allowing the adjacent circles to separate.
Suppose $M$ is a
surface with a triangulation $\mathcal{T}=\{V, E, F\}$,
where $V, E, F$ are the sets of vertices, edges and faces respectively.
We use $i$, $\{ij\}$, $\{ijk\}$ to denote a vertex, an edge and a face respectively, where
$i,j,k$ are natural numbers.
A weight on the triangulated surface is a map $I: E\rightarrow (-1, +\infty)$. We use $I_{ij}$ to denote $I(\{ij\})$ for simplicity.
A weighted triangulated surface is
denoted by $(M, \mathcal{T}, I)$ in this paper.

\begin{definition}\label{definition of inversive distance circle packing}
Suppose $(M, \mathcal{T}, I)$ is a weighted triangulated surface.
An inversive distance circle packing metric is a map $r: V\rightarrow (0, +\infty)$ such that
\begin{description}
  \item[(1)] The edge length $l_{ij}$ of $\{ij\}$ is
  \begin{equation}\label{Euclidean length}
  l_{ij}=\sqrt{r_i^2+r_j^2+2r_ir_jI_{ij}}
  \end{equation}
  for Euclidean background geometry and
   \begin{equation}\label{hyperbolic length}
  l_{ij}=\cosh^{-1}(\cosh r_i\cosh r_j+I_{ij}\sinh r_i\sinh r_j)
  \end{equation}
  for hyperbolic background geometry;
  \item[(2)] With the assignment of edge lengths $l_{ij}, l_{jk}, l_{ik}$
  by $(\ref{Euclidean length})$ $($respectively $(\ref{hyperbolic length})$$)$,
  the triangle $\{ijk\}$ could be embedded in $2$-dimensional Euclidean space $\mathbb{E}^2$
  (respectively $2$-dimensional hyperbolic space $\mathbb{H}^2$) as a nondegenerate triangle.
\end{description}
\end{definition}

The condition (2) in Definition \ref{definition of inversive distance circle packing} is called a nondegenerate condition.
If two circles $C_i$ and $C_j$ with radii $r_i$ and $r_j$ respectively
are put in the plane with $l_{ij}$ as the distance of the centers of $C_i, C_j$, then the inversive distance of the two circles is $I_{ij}$.
If $I_{ij}\in [0,1]$ for any edge $\{ij\}\in E$, the inversive distance circle packing is reduced to Thurston's circle packing \cite{T1}.
If $I_{ij}\in (-1, 1]$ for any edge $\{ij\}\in E$, the inversive distance circle packing is reduced to the circle packing studied by Zhou \cite{Z}.
If $I_{ij}\in [0,+\infty)$ for any edge $\{ij\}\in E$, the inversive distance circle packing was studied by Guo \cite{Guo} and Luo \cite{L}.
For more information on inversive distance circle packing, see \cite{BH,BS,Guo,S}.

Bowers-Stephenson \cite{BS} conjectured that the inversive distance circle packing on surfaces is rigid.
The infinitesimal rigidity and global rigidity were proved
by Guo \cite{Guo} and Luo \cite{L} respectively for circle packings with nonnegative inversive distance,
which generalized the Andreev-Thurston rigidity theorem.
Following the proof in \cite{Guo,L}, the author \cite{X2} proved the rigidity of circle packing
for inversive distance in $(-1, +\infty)$ recently.

\begin{theorem}[\cite{Guo,L,X2}]\label{main theorem}
Suppose $(M, \mathcal{T}, I)$ is a weighted triangulated surface
with the weight $I: E\rightarrow (-1, +\infty)$ satisfying the structure condition
\begin{equation}\label{structure condition}
\begin{aligned}
I_{ij}+I_{ik}I_{jk}\geq 0, I_{ik}+I_{ij}I_{jk}\geq 0, I_{jk}+I_{ij}I_{ik}\geq 0,\ \ \forall\{ijk\}\in F.
\end{aligned}
\end{equation}
Then the inversive distance circle packing metric on $(M, \mathcal{T}, I)$ is
uniquely determined by the discrete curvature (up to scaling for the Euclidean background geometry).
\end{theorem}

The basic strategy in \cite{Guo,L,X2} to prove Theorem  \ref{main theorem} is
to apply the variational principle introduced by de Verdi\`{e}re \cite{DV} to inversive distance circle packing,
which could be separated into the following three steps.
The first step is to prove the admissible space of inversive distance packing metrics for a single
triangle is simply connected; The second step is to prove that the Jacobian matrix of the inner angles of a triangle in terms of
some appropriate parametrization of the circle radii
is symmetric and negative semi-definite (or negative definite), which ensures the definition of a locally concave function;
The third step is to extend the locally concave function to be a globally defined concave function,
which has been systematically studied in \cite{BPS,L}, and use this concave function to prove the rigidity.

In this paper, we give a new proof of Theorem \ref{main theorem}.
In the first step, the proof in \cite{Guo,L,X2} for simply connectivity of the admissible space for a single triangle
is based on the triangle inequalities, which can not be
generalized to three or higher dimensional cases.
In this paper, we give a new proof of the simply connectivity using the Cayley-Menger determinant, which could be used to
characterize the nondegeneracy of a simplex in any dimension.
This proof enables us to prove the simply connectivity
of admissible space of Thurston's sphere packing metrics for a single tetrahedron in three dimension \cite{HX1,HX2}.
In the second step,
the arguments in \cite{Guo,X2} to prove the negative semi-definiteness (or negative definiteness) of
the Jacobian matrix of inner angles in a triangle
is based on a lengthy estimate of the eigenvalues of the matrix under the nondegenerate condition.
In this paper, we give a new and short proof of the negative definiteness involving only the rank of the Jacobian matrix and
connectivity of the parameterized admissible space for a triangle,
which greatly simplifies the arguments in \cite{Guo,X2}. The third step is the same as that in \cite{Guo,L,X2}.

In this paper, we only study the rigidity of inversive distance circle packing in Euclidean and
hyperbolic background geometry. For the rigidity
of inversive distance circle packing in spherical background geometry, see \cite{BB,BBP,MS}.
Deformation of inversive distance circle packing metrics
on surfaces by discrete curvature flows was also studied recently, see \cite{GJ1,GJ2,GJ3,GX1}.
Inversive distance circle packing has lots of practical applications, see \cite{BH,HS,ZG,ZZGLG}.

This paper is organized as follows.
In Section \ref{section 2}, we give a new proof of Theorem \ref{main theorem} in Euclidean background geometry.
In Section \ref{section 3}, we give a new proof of Theorem \ref{main theorem} in hyperbolic background geometry.

\section{Rigidity of Euclidean inversive distance circle packing}\label{section 2}
\subsection{Admissible space of Euclidean inversive distance circle packing metrics for a single triangle}

Suppose $\sigma=\{123\}\in F$ is a topological triangle in $(M, \mathcal{T},I)$.
The corresponding edge set of the triangle is denoted by $E_{\sigma}=\{\{12\}, \{13\}, \{23\}\}$.
We denote $\eta$ as the restriction of
the weight $I: E \rightarrow (-1, +\infty)$ on the edge set $E_{\sigma}$.
Given a weight $\eta$ on the edge set $E_\sigma$ satisfying the structure condition (\ref{structure condition}),
the admissible space $\Omega^E_{123}(\eta)$ of Euclidean inversive distance circle packing metrics for the triangle $\{123\}$
is defined to be the set of Euclidean inversive distance circle packing metrics $(r_1, r_2, r_3)\in \mathbb{R}^3_{>0}$ such that
the triangle $\{123\}$ with the edge lengths given by (\ref{Euclidean length})
exists in 2-dimensional Euclidean space $\mathbb{E}^2$.

To simplify the notations, we set
\begin{equation*}
I_i=I_{jk},  \ \{i,j,k\}=\{1,2,3\}.
\end{equation*}
Then $l_{ij}^2=r_i^2+r_j^2+2r_ir_jI_k$.
Set
\begin{equation*}
\begin{aligned}
G_0(l)=\left(
         \begin{array}{cccc}
           0 & 1 & 1 & 1 \\
           1 & 0 & l_{12}^2 & l_{13}^2 \\
           1 & l_{12}^2 & 0 & l_{23}^2 \\
           1 & l_{13}^2 & l_{23}^2 & 0 \\
         \end{array}
       \right)
\end{aligned}
\end{equation*}
to be the Cayley-Menger $4\times 4$-matrix.
Recall the following result characterizing
the nondegeneracy of a Euclidean triangle $\{123\}$ with positive edge lengths $l_{12}, l_{13}, l_{23}$.

\begin{lemma}[\cite{TW}, Proposition 2.4.1]
A triangle with positive edge lengths $l_{12}, l_{13}, l_{23}$ exists in $\mathbb{E}^2$ if and only
if $\det G_0(l)<0$.
\end{lemma}
\begin{remark}
By direct calculations, we have
\begin{equation*}
\begin{aligned}
\det G_0(l)=-(l_{12}+l_{13}+l_{23})(l_{12}+l_{13}-l_{23})(l_{12}+l_{23}-l_{13})(l_{13}+l_{23}-l_{12}),
\end{aligned}
\end{equation*}
which implies $\det G_0(l)<0$ is equivalent to the triangle inequalities. This was also observed in \cite{Guo}.
The advantage of using $\det G_0(l)<0$ to characterize the nondegeneracy is that
we just need one inequality $\det G_0(l)<0$ instead of three triangle inequalities.
Furthermore, this characterization of nondegeneracy of simplex could be generalized
to high dimensional case \cite{TW}.
\end{remark}

Submitting (\ref{Euclidean length}) into $\det G_0(l)$, we have
\begin{equation*}
\begin{aligned}
\det G_0(l)=&-4[r_1^2r_2^2(1-I^2_{3})+r_1^2r_3^2(1-I^2_{2})+r_2^2r_3^2(1-I^2_{1})\\
 &+2r_1^2r_2r_3(I_{1}+I_{2}I_{3})+2r_1r_2^2r_2(I_{2}+I_{1}I_{3})
 +2r_1r_2r_3^2(I_{3}+I_{1}I_{2})].
\end{aligned}
\end{equation*}
Set
\begin{equation*}
\gamma_{i}=I_i+I_jI_k, \ \kappa_i=r_i^{-1}
\end{equation*}
and
\begin{equation*}
Q=\kappa_1^2(1-I^2_{1})+\kappa_1^2(1-I^2_{2})+\kappa_3^2(1-I^2_{3})
 +2\kappa_1\kappa_2\gamma_{3}+2\kappa_1\kappa_3\gamma_{2}+2\kappa_2\kappa_3\gamma_{1},
 \end{equation*}
then we have
\begin{equation*}
\begin{aligned}
\det G_0(l)=-4r_1^2r_2^2r_3^2Q.
\end{aligned}
\end{equation*}
\begin{lemma}[\cite{Guo,X2,Z}]\label{nondegenerate lemma Euclidean}
A Euclidean triangle $\{123\}$ with edge lengths $l_{12}, l_{13}, l_{23}$
given by $(\ref{Euclidean length})$ exists in $\mathbb{E}^2$ if and only if $Q>0$.
\end{lemma}
Set
\begin{equation}\label{hi}
\begin{aligned}
h_1=&\kappa_1(1-I_1^2)+\kappa_2\gamma_{3}+\kappa_3\gamma_{2},\\
h_2=&\kappa_2(1-I_2^2)+\kappa_1\gamma_{3}+\kappa_3\gamma_{1},\\
h_3=&\kappa_3(1-I_3^2)+\kappa_1\gamma_{2}+\kappa_2\gamma_{1}.
\end{aligned}
\end{equation}
By Lemma \ref{nondegenerate lemma Euclidean}, $(r_1,r_2,r_3)\in \mathbb{R}^3_{>0}$ generates a degenerate Euclidean triangle if and only if
\begin{equation}\label{degenerate condition}
\begin{aligned}
Q=\kappa_1h_1+\kappa_2h_2+\kappa_3h_3\leq 0.
\end{aligned}
\end{equation}
\begin{remark}\label{relation of hi and h-jk,i}
For a nondegenerate inversive distance circle packing metric of the triangle $\{123\}$,
there exists a geometric center $C_{123}$ of the triangle $\{123\}$, which has the same circle power
to the circles attached to the vertices $\{1,2,3\}$.
$h_i$ in $(\ref{hi})$ is a positive multiplication of the signed distance
$h_{jk,i}$ of $C_{123}$ to the edge $\{jk\}$, which is defined to be positive
if $C_{123}$ is on the same side of the line determined by $\{jk\}$ as the triangle $\{123\}$  and negative otherwise
(or zero if $C_{123}$ is on the line).
By direct calculations, we have
\begin{equation}\label{calculation of h_ij,k}
\begin{aligned}
h_{ij,k}
=\frac{r_1^2r_2^2r_3^2}{2l_{ij}A_{123}}[\kappa_k^2(1-I_k^2)+\kappa_j\kappa_k\gamma_{i}+\kappa_i\kappa_k\gamma_{j}]
=\frac{r_1^2r_2^2r_3^2}{2l_{ij}A_{123}}\kappa_kh_k,
\end{aligned}
\end{equation}
where $A_{123}$ is the area of the triangle $\{123\}$.
Note that $h_1,h_2,h_3$ are well-defined for any $(r_1, r_2, r_3)\in \mathbb{R}^3_{>0}$, while $h_{12,3}, h_{13,2}, h_{23,1}$
are defined for nondegenerate inversive distance circle packing metrics.
Refer to \cite{G1,G2,GT,T} for more information on the geometric center of triangles.
\end{remark}

Suppose $(r_1, r_2, r_3)\in \mathbb{R}^3_{>0}$ is a degenerate inversive distance circle packing metric,
then one of the following two cases happens by (\ref{degenerate condition}).
\begin{description}
  \item[(1)] At least one of $h_1,h_2,h_3$ is zero;
  \item[(2)] None of $h_1,h_2,h_3$ is zero.
\end{description}
We will prove that case (1) never happens. Furthermore,
we will prove that only one of $h_i,h_j,h_k$ is negative and the others are positive in case (2).

Note that $Q\leq 0$ is equivalent to the following quadratic inequality of $\kappa_i$
\begin{equation}
\begin{aligned}
A_i\kappa_i^2+B_i\kappa_i+C_i\geq 0,
\end{aligned}
\end{equation}
where
\begin{equation}\label{A1,B1,C1}
\begin{aligned}
A_i=&I^2_{i}-1,\\
B_i=&-2(\kappa_j\gamma_k+\kappa_k\gamma_j)\leq0,\\
C_i=&\kappa_j^2(I_j^2-1)+\kappa_k^2(I^2_{k}-1)-2\kappa_j\kappa_k\gamma_i.
\end{aligned}
\end{equation}
with $\{i,j,k\}= \{1,2,3\}$.
By direct calculations, the determinant $\Delta_i=B_i^2-4A_iC_i$ is given by
\begin{equation}\label{Delta1 Euclidean}
\begin{aligned}
\Delta_i
=4(I_1^2+I_2^2+I_3^2+2I_1I_2I_3-1)(\kappa_j^2+\kappa_k^2+2\kappa_j\kappa_kI_i).
\end{aligned}
\end{equation}

\begin{lemma}\label{positive Delta Euclidean}
Suppose $\eta=(I_1, I_2, I_3)$ is a weight on the edges of a triangle $\{123\}$
satisfying the structure condition $(\ref{structure condition})$.
If $I_i>1$, then $\Delta_i>0$.
\end{lemma}

\proof
It is straight forward to check that $\kappa_j^2+\kappa_k^2+2\kappa_j\kappa_kI_i>0$.
We just need to check $I_1^2+I_2^2+I_3^2+2I_1I_2I_3-1>0$ by (\ref{Delta1 Euclidean}).
If $I_j\geq 0, I_k\geq 0$, then
$$I_1^2+I_2^2+I_3^2+2I_1I_2I_3-1>I_j^2+I_k^2+2I_1I_2I_3\geq 0.$$
If $I_j<0$, then $I_j\in (-1,0)$ and
$$I_1^2+I_2^2+I_3^2+2I_1I_2I_3-1=(I_k+I_iI_j)^2+(1-I_j^2)(I_i^2-1)>0.$$
Similar argument applies for the case $I_k<0$.
Therefore, under the structure condition (\ref{structure condition}) and $I_i>1$,
we have $I_1^2+I_2^2+I_3^2+2I_1I_2I_3-1>0$. \qed

Now we can prove that the case (1) never happens.

\begin{lemma}\label{lemma none of h1,h2,h3 is 0}
Suppose $(r_1, r_2, r_3)\in \mathbb{R}^3_{>0}$ is a degenerate Euclidean inversive distance circle packing
metric for a triangle $\{123\}$ with a weigh $\eta: E_\sigma\rightarrow (-1, +\infty)$
satisfying the structure condition $(\ref{structure condition})$,
then none of $h_1, h_2, h_3$ is zero.
\end{lemma}
\proof
We prove the lemma by contradiction.
By the degenerate condition (\ref{degenerate condition}),
if one of $h_1, h_2, h_3$ is zero, then there is another one of $h_1, h_2, h_3$ that is nonpositive.
Without loss of generality, we assume $h_1=0, h_2\leq 0$ for
a degenerate Euclidean inversive distance circle packing metric $(r_1, r_2, r_3)\in \mathbb{R}^3_{>0}$.

By $h_1=0$, we have
$\kappa_1(I_1^2-1)=\kappa_2\gamma_3+\kappa_3\gamma_2,$
which implies $I_1\geq 1$ by the structure condition (\ref{structure condition}).
If $I_1>1$, then we can take $Q\leq 0$ as a quadratic inequality in $\kappa_1$
\begin{equation}\label{Q quadratic nonpositive in proof Euclidean}
\begin{aligned}
A_1\kappa_1^2+B_1\kappa_1+C_1\geq 0,
\end{aligned}
\end{equation}
where $A_1, B_1, C_1$ are given by (\ref{A1,B1,C1}) with $A_1=I_1^2-1>0$.
By Lemma \ref{positive Delta Euclidean}, we have $\Delta_1>0$. Then (\ref{Q quadratic nonpositive in proof Euclidean}) implies
\begin{equation*}
\begin{aligned}
\kappa_1\geq \frac{-B_1+\sqrt{\Delta_1}}{2A_1}  \ \ \ \mbox{or} \ \ \ \kappa_1\leq \frac{-B_1-\sqrt{\Delta_1}}{2A_1},
\end{aligned}
\end{equation*}
which is equivalent to
\begin{equation*}
\begin{aligned}
-2h_1=2A_1\kappa_1+B_1\geq \sqrt{\Delta_1}>0  \ \ \ \mbox{or} \ \ \ -2h_1=2A_1\kappa_1+B_1\leq -\sqrt{\Delta_1}<0.
\end{aligned}
\end{equation*}
This contradicts $h_1=0$. Therefore, $I_1=1$.
By $h_1=0$ again, we have $\gamma_2=\gamma_3=0$, which implies $I_2+I_3=0$.

By $h_2=\kappa_2(1-I_2^2)+\kappa_1\gamma_3+\kappa_3\gamma_1\leq0$, we have $I_2\geq1$,
which implies $I_3=-I_2\leq-1$. This is impossible. \qed

By Lemma \ref{lemma none of h1,h2,h3 is 0},
if $(r_1, r_2, r_3)\in \mathbb{R}^3_{>0}$ is a degenerate Euclidean inversive distance circle packing
metric for a triangle $\{123\}$, at least one of $h_1,h_2,h_3$ is negative and the others are nonzero.
Furthermore, we have the following result.

\begin{lemma}\label{lemma no h1<0 h2<0 Euclidean}
Suppose $\{123\}$ is
a triangle with a weigh $\eta: E_\sigma\rightarrow (-1, +\infty)$
satisfying the structure condition $(\ref{structure condition})$ and
$(r_1, r_2, r_3)\in \mathbb{R}^3_{>0}$. Then there exists no subset $\{i,j\}\subset \{1,2,3\}$ such that
$h_i<0$ and $h_j<0$.
\end{lemma}
\proof
Without loss of generality, we consider the case $h_1<0, h_2<0$.
By $h_1<0, h_2<0$, we have
\begin{equation*}
\begin{aligned}
(I_1^2-1)\kappa_1>\kappa_2\gamma_3+\kappa_3\gamma_2,\ \
(I_2^2-1)\kappa_2>\kappa_1\gamma_3+\kappa_3\gamma_1,
\end{aligned}
\end{equation*}
which implies $I_1>1, I_2>1$ and $(I_1^2-1)(I_2^2-1)>\gamma_3^2$ by the structure condition (\ref{structure condition}).
Note that
\begin{equation*}
\begin{aligned}
(I_1^2-1)(I_2^2-1)-\gamma_3^2
=-I_1^2-I_2^2-I_3^2-2I_1I_2I_3+1<0
\end{aligned}
\end{equation*}
by the proof of Lemma \ref{positive Delta Euclidean}. This is a contradiction.\qed

\begin{remark}
No matter $(r_1, r_2, r_3)\in \mathbb{R}^3_{>0}$ is a nondegenerate or degenerate inversive distance
circle packing metric for the triangle $\{123\}$, Lemma $\ref{lemma no h1<0 h2<0 Euclidean}$ is valid. For nondegenerate
inversive distance circle packing metrics, Lemma $\ref{lemma no h1<0 h2<0 Euclidean}$ implies that
the geometric center does not lie in some special regions in the plane relative to the triangle.
\end{remark}

Now we can prove the main result of this subsection.

\begin{proposition}[\cite{Guo,X2}]\label{simply connect of admi space with weight}
Suppose $\sigma=\{123\}\in F$ is a triangle in $(M, \mathcal{T})$ with a weight $\eta: E_\sigma\rightarrow (-1, +\infty)$
satisfying the structure condition $(\ref{structure condition})$.  Then the admissible
space $\Omega^E_{123}(\eta)$ of Euclidean inversive distance circle packing metrics $(r_1, r_2, r_3)\in \mathbb{R}^3_{>0}$
is nonempty and simply connected.
Furthermore, the set of degenerate inversive distance circle packing metric is a disjoint union
$\cup_{i\in P}V_i,$
where $P=\{i\in \{1,2,3\}|I_i>1\}$ and
$$V_i=\left\{(r_1, r_2, r_3)\in \mathbb{R}^3_{>0}|\kappa_i\geq \frac{-B_i+\sqrt{\Delta_i}}{2A_i}\right\}$$
is bounded by an analytic graph on $\mathbb{R}^2_{>0}$.
\end{proposition}
\proof
Suppose $(r_1, r_2, r_3)\in \mathbb{R}^3_{>0}$ is a degenerate inversive distance circle packing metric
for the triangle $\{123\}$, then we have
$Q=\kappa_1h_1+\kappa_2h_2+\kappa_3h_3\leq 0.$
By Lemma \ref{lemma none of h1,h2,h3 is 0} and Lemma \ref{lemma no h1<0 h2<0 Euclidean},
one of $h_1,h_2,h_3$ is negative and the others are positive.

Suppose $h_i<0$ and $h_j>0, h_k>0$ with $\{i,j,k\}=\{1,2,3\}$.
Then we have $I_i>1$ by $h_i<0$.
Take $Q\leq 0$ as a quadratic inequality of $\kappa_i$, we have
$A_i\kappa_i^2+B_i\kappa_i+C_i\geq 0,$
where $A_i, B_i, C_i$ are defined by (\ref{A1,B1,C1}) with $A_i=I_i^2-1>0$.
By Lemma \ref{positive Delta Euclidean}, we have $\Delta_i>0$. Then $A_i\kappa_i^2+B_i\kappa_i+C_i\geq 0$ implies
\begin{equation*}
\begin{aligned}
\kappa_i\geq \frac{-B_i+\sqrt{\Delta_i}}{2A_i}  \ \ \ \mbox{or} \ \ \ \kappa_i\leq \frac{-B_i-\sqrt{\Delta_i}}{2A_i}.
\end{aligned}
\end{equation*}
Note that $h_i<0$ is equivalent to $\kappa_i>\frac{-B_i}{2A_i}$.
This implies $\kappa_i\geq \frac{-B_i+\sqrt{\Delta_i}}{2A_i}$. Therefore, the set of degenerate
inversive distance circle packing metrics is contained in $\cup_{i\in P}V_i$.

On the other hand, if $I_i>1$, then for any $(r_1, r_2, r_3)\in V_i$, we have $Q\leq 0$, which implies
any element $(r_1, r_2, r_3)\in V_i$ is a degenerate inversive distance circle packing metric.
Therefore,
$\Omega^E_{123}(\eta)=\mathbb{R}^3_{>0}\setminus \cup_{i\in P}V_i,$
where $P=\{i\in \{1,2,3\}|I_i>1\}$.

For any $(r_1, r_2, r_3)\in V_i$, we have $\kappa_i>\frac{-B_i}{2A_i}$, which is equivalent to $h_i<0$.
This implies $V_i\cap V_j=\emptyset$ if $I_i>1$ and $I_j>1$ by Lemma \ref{lemma no h1<0 h2<0 Euclidean}.

Note that $V_i$ is bounded by an analytic graph on $\mathbb{R}^2_{>0}$
$$V_i=\{(r_1, r_2, r_3)\in \mathbb{R}^3_{>0}|r_i\leq \frac{2A_i}{-B_i+\sqrt{\Delta_i}}\}.$$
This implies $\Omega^E_{123}(\eta)=\mathbb{R}^3_{>0}\setminus \cup_{i\in P}V_i$ is homotopy equivalent to $\mathbb{R}^3_{>0}$.
Therefore, $\Omega^E_{123}(\eta)$ is simply connected. \qed

\begin{remark}\label{Vi hi<0, hj,hk>0 Euclidean}
Suppose $\sigma=\{123\}\in F$ is a triangle with a weight $\eta: E_\sigma\rightarrow (-1, +\infty)$
satisfying the structure condition $(\ref{structure condition})$.
For $(r_1, r_2, r_3)\in V_i$, we have $h_i<0$ and $h_j>0, h_k>0$.
\end{remark}

\begin{remark}
The simply connectivity of the admissible space of nondegenerate inversive distance circle packing
metrics was first proved by Guo \cite{Guo} for nonnegative inversive distance and then by the author \cite{X2}
for inversive distance in $(-1, +\infty)$ satisfying the structure condition (\ref{structure condition}).
The proof of simply connectivity presented here
is motivated by the proof of simply connectivity of admissible space of
sphere packing metrics of a tetrahedron in 3-dimension \cite{GH,GJS, X1}.
The advantage of the proof of Proposition \ref{simply connect of admi space with weight} is that
we have a precise description of the boundary of the admissible space,
each connected component of which is an analytic graph on $\mathbb{R}^2_{>0}$, and
the proof could be generalized to 3-dimensional case to prove the
simply connectivity of admissible space of Thurston's sphere packing metrics for a tetrahedron \cite{GH,GJS,HX1,HX2,X1}.
\end{remark}

In a single triangle $\{123\}$, we denote $\theta_i$ as the angle at vertex $i$.
We have the following result.
\begin{lemma}[\cite{L,X2}]\label{continuity of Euclidan inner angle}
Suppose $\{123\}\in F$ is a triangle
with a weight $\eta: E_\sigma\rightarrow (-1, +\infty)$
satisfying the structure condition $(\ref{structure condition})$.
Then $\theta_1, \theta_2, \theta_3$ defined for $(r_1,r_2,r_3)\in\Omega^E_{123}(\eta)$
could be extended by constants to be continuous functions
$\widetilde{\theta}_1, \widetilde{\theta}_2, \widetilde{\theta}_3$ defined on $\mathbb{R}^3_{>0}$.
\end{lemma}
\proof
If $r=(r_1,r_2,r_3)\in \Omega^E_{123}(\eta)$ tends to a point $\overline{r}=(\overline{r}_1,\overline{r}_2,\overline{r}_3)$
in the boundary $\partial V_i$ of $V_i$ in $\mathbb{R}^3_{>0}$,
we have the area $A_{123}\rightarrow 0$, $l_{ij}(r)\rightarrow l_{ij}(\overline{r})>0$ and
$l_{ik}(r)\rightarrow l_{ik}(\overline{r})>0$,
which implies $\sin\theta_i=\frac{2A_{123}}{l_{ij}l_{ik}}\rightarrow 0$. Therefore, $\theta_i\rightarrow \pi$ or $0$.

Take $u_i=\ln r_i$.
By Lemma \ref{symmetry Euclidean}, we have
$$\frac{\partial \theta_i}{\partial u_i}=-\frac{\partial \theta_j}{\partial u_i}-\frac{\partial \theta_k}{\partial u_i}
=-\frac{h_{ij,k}}{l_{ij}}-\frac{h_{ik,j}}{l_{ik}}.$$
As $h_j>0, h_k>0$ for $\overline{r}\in \partial V_i$ by Remark \ref{Vi hi<0, hj,hk>0 Euclidean},
we have $\frac{\partial \theta_i}{\partial u_i}<0$ for $r\in \Omega^E_{123}(\eta)$ around $\overline{r}$ by (\ref{calculation of h_ij,k}).
Therefore, if we increase $r_i$ at $\overline{r}$, which results the triangle does not degenerate, we shall have $\theta_i$
decrease. This implies $\theta_i\rightarrow \pi$ as $(r_1,r_2,r_3)\rightarrow (\overline{r}_1,\overline{r}_2,\overline{r}_3)$.
By $\theta_i+\theta_j+\theta_k=\pi$, we have $\theta_j, \theta_k\rightarrow 0$.

Then we can extend $\theta_1, \theta_2,\theta_3$ defined on $\Omega^E_{123}(\eta)$ to be continuous functions defined on $\mathbb{R}^3_{>0}$
by setting
\begin{equation*}
\begin{aligned}
\widetilde{\theta}_i(r_1,r_2,r_3)=\left\{
                       \begin{array}{ll}
                         \theta_i, & \hbox{if $(r_1,r_2,r_3)\in \Omega^E_{123}(\eta)$;} \\
                         \pi, & \hbox{if $(r_1,r_2,r_3)\in V_i$;} \\
                         0, & \hbox{otherwise.}
                       \end{array}
                     \right.
\end{aligned}
\end{equation*}
\qed

Denote
\begin{equation*}
\begin{aligned}
\Gamma=\{(I_1, I_2, I_3)\in (-1, +\infty)^3|\gamma_1\geq 0,\gamma_2\geq 0,\gamma_3\geq 0\}
\end{aligned}
\end{equation*}
as the space of weights on the edges of a triangle $\{123\}$ satisfying the structure condition (\ref{structure condition}).

\begin{lemma}\label{weight space is connected}
$\Gamma$ is connected.
\end{lemma}
\proof
It is obviously that $[0, +\infty)^3\subset \Gamma$.
By $\gamma_i=I_i+I_jI_k\geq 0, i=1,2,3$, we have $I_1+I_2\geq 0, I_1+I_3\geq 0, I_2+I_3\geq 0$,
which implies that at most one of $I_1, I_2, I_3$ is negative.
Without loss of generality, we consider the case $(I_1, I_2, I_3)\in \Gamma$ with $I_1<0, I_2\geq 0, I_3\geq 0$.
It is straight forward to check that
$(tI_1, I_2, I_3)\in \Gamma$ for any $t\in [0,1]$, which implies $\Gamma$ is connected. \qed

Using the space $\Gamma$,
we can further define the following 6-dimensional parameterized admissible space
\begin{equation*}
\begin{aligned}
\Omega^E_{123}=\cup_{\eta\in \Gamma}\Omega^E_{123}(\eta).
\end{aligned}
\end{equation*}

\begin{lemma}\label{connectivity of Omega123 Euclidean}
$\Omega^E_{123}$ is connected.
\end{lemma}
\proof
Suppose $\eta_0\in \Gamma$, then there exists $r_0\in \Omega^E_{123}(\eta_0)$
with $Q(\eta_0, r_0)>0$ by the nonempty property of $\Omega^E_{123}(\eta_0)$ in Proposition \ref{simply connect of admi space with weight}.
Consider the continuous function $Q(\eta, r_0)$ of $\eta$. As $Q(\eta_0, r_0)>0$,
there is a connected neighborhood $U_{\eta_0}\subset \Gamma$ of $\eta_0$ such
that $Q(\eta, r_0)>0$ for any $\eta\in U_{\eta_0}$.
This implies that for any $\eta\in U_{\eta_0}$,
any two points $(\eta, r_A)\in\Omega^E_{123}$ and $(\eta_0, r_B)\in\Omega^E_{123}$
could be connected by a path in $\Omega^E_{123}$ by Proposition \ref{simply connect of admi space with weight}.
In this case, we call the space $\Omega^E_{123}(\eta)$ and $\Omega^E_{123}(\eta_0)$ could be connected by a path in $U_{\eta_0}$.
Taking $\Omega^E_{123}(\eta)$ as a point.
Then for any $\eta_A, \eta_B\in \Gamma$, the existence of a path from $\Omega^E_{123}(\eta_A)$ to $\Omega^E_{123}(\eta_B)$ in  $\Gamma$ follows from the connectivity of $\Gamma$ and finite covering theorem, which implies that $\Omega^E_{123}$ is connected. \qed

\subsection{Negative semi-definiteness of the Jacobian matrix}

Set $u_i=\ln r_i$.
The following result on the matrix $\Lambda^E_{123}=\frac{\partial(\theta_1, \theta_2, \theta_3)}{\partial (u_1, u_2, u_3)}$
is known.

\begin{lemma}[\cite{CL,Guo,T1,X2,Z}]\label{symmetry Euclidean}
Suppose $(r_1,r_2,r_3)$ is a nondegenerate Euclidean inversive distance circle packing metric
for a triangle $\{123\}$ with a weight $\eta$ satisfying the structure condition $(\ref{structure condition})$, then
\begin{equation}\label{derivative of theta}
\frac{\partial \theta_i}{\partial u_j}=\frac{\partial \theta_j}{\partial u_i}
=\frac{r_1^2r_2^2r_3^2}{2A_{123}l_{ij}^2}\left[\kappa_k^2(1-I_k^2)+\kappa_j\kappa_k\gamma_{i}+\kappa_i\kappa_k\gamma_{j}\right]
\end{equation}
for any adjacent vertices $i, j$, where $A_{123}$ is the area of the triangle $\{123\}$, and
\begin{equation}\label{partial theta i of u i}
\frac{\partial \theta_i}{\partial u_i}=-\frac{\partial \theta_i}{\partial u_j}-\frac{\partial \theta_i}{\partial u_k}.
\end{equation}
Specially, for $\eta=(I_1, I_2, I_3)=(1,1,1)\in \Gamma$,
the Jacobian matrix $\Lambda^E_{123}=\frac{\partial(\theta_1, \theta_2, \theta_3)}{\partial (u_1, u_2, u_3)}$
is negative semi-definite with a zero eigenvalue and two negative eigenvalues.
\end{lemma}

\begin{remark}
By $(\ref{calculation of h_ij,k})$, $(\ref{partial theta i of u i})$ and Remark \ref{Vi hi<0, hj,hk>0 Euclidean},
if $\eta\in \Gamma$ and $(r_1,r_2,r_3)\in \Omega^E_{123}(\eta)$
tends to a point $(\overline{r}_1, \overline{r}_2, \overline{r}_3)\in \partial V_i$,
we have $\frac{\partial \theta_i}{\partial u_j}\rightarrow +\infty$, $\frac{\partial \theta_i}{\partial u_k}\rightarrow +\infty$ and
$\frac{\partial \theta_i}{\partial u_i}\rightarrow -\infty$.
\end{remark}

Set $d_{ij}=\frac{\partial l_{ij}}{\partial u_i}$ and $d_{ji}=\frac{\partial l_{ij}}{\partial u_j}$.
Then $d_{ij}=\frac{r_i(r_i+r_jI_k)}{l_{ij}}$, $d_{ji}=\frac{r_j(r_j+r_iI_k)}{l_{ij}}$ and $d_{ij}+d_{ji}=l_{ij}$.
This is a type of conformal metric studied in \cite{G1,G2,GT,T}.

\begin{lemma}[\cite{Guo,X2}]\label{negativity of Lambda}
Suppose $r=(r_1, r_2, r_3)\in \mathbb{R}^3_{>0}$ is a nondegenerate Euclidean inversive distance circle packing metric
for a triangle $\{123\}$ with a weight $\eta$ satisfying the structure condition $(\ref{structure condition})$, then
the Jacobian matrix $\Lambda^E_{123}(\eta)=\frac{\partial(\theta_1, \theta_2, \theta_3)}{\partial (u_1, u_2, u_3)}$
is negative semi-definite with one dimensional kernel $\{t(1,1,1)|t\in \mathbb{R}\}$.
\end{lemma}
\proof
By the chain rules, we have
\begin{equation*}
\begin{aligned}
\Lambda^E_{123}
=\frac{\partial(\theta_1, \theta_2, \theta_3)}{\partial (l_{23}, l_{13}, l_{12})}
\cdot
\frac{\partial (l_{23},l_{13},l_{12})}{\partial (u_1, u_2, u_3)}.
\end{aligned}
\end{equation*}
By direct calculations,
\begin{equation*}
\begin{aligned}
\det \frac{\partial (l_{23},l_{13},l_{12})}{\partial (u_1, u_2, u_3)}
=&\det \left(
          \begin{array}{ccc}
           0 & d_{23} & d_{32} \\
            d_{13} & 0 & d_{31} \\
            d_{12} & d_{21} & 0 \\
          \end{array}
        \right)\\
=&\frac{r_1r_2r_3}{l_{12}l_{13}l_{23}}[2r_1r_2r_3(1+I_1I_2I_3)+r_1(r_2^2+r_3^2)(I_1+I_2I_3)\\
&+r_2(r_1^2+r_3^2)(I_2+I_1I_3)+r_3(r_1^2+r_2^2)(I_3+I_1I_2)]\\
\geq&\frac{2r_1^2r_2^2r_3^2}{l_{12}l_{13}l_{23}}(1+I_1I_2I_3+I_1+I_2I_3+I_2+I_1I_3+I_3+I_1I_2)\\
=&\frac{2r_1^2r_2^2r_3^2}{l_{12}l_{13}l_{23}}(1+I_1)(1+I_2)(1+I_3)>0,
\end{aligned}
\end{equation*}
which implies the matrix $\frac{\partial (l_{23},l_{13},l_{12})}{\partial (u_1, u_2, u_3)}$
is nondegenerate. Therefore,  the rank of $\Lambda^E_{123}$ is the same as that of the matrix
$\frac{\partial(\theta_1, \theta_2, \theta_3)}{\partial (l_{23}, l_{13}, l_{12})}$,
which is known  to be 2 for nondegenerate Euclidean triangles.

Taking $\Lambda^E_{123}$ as a matrix-valued function defined on $\Omega^E_{123}$, then
the two nonzero eigenvalues of $\Lambda^E_{123}$ are continuous functions of $(\eta, r)\in \Omega^E_{123}$.
By the connectivity of $\Omega^E_{123}$ in Lemma \ref{connectivity of Omega123 Euclidean},
the nonzero eigenvalues do not change sign in $\Omega^E_{123}$.
Note that $\Lambda^E_{123}(\eta_0)$ is negative semi-definite with two negative eigenvalues for $\eta_0=(1,1,1)$ by Lemma \ref{symmetry Euclidean}.
This implies that
$\Lambda^E_{123}(\eta)$ is negative semi-definite with two negative eigenvalues for any $\eta\in \Gamma$.
The kernel of $\Lambda^E_{123}(\eta)$ is $\{t(1,1,1)|t\in \mathbb{R}\}$ follows from the scaling invariance of $\theta_1, \theta_2, \theta_3$.
\qed

\begin{remark}
The negative semi-definiteness of the Jacobian matrix
$\Lambda^E_{123}$ for Thurston's circle packing metric is well-known, see \cite{CL, T1, Z}.
The negative semi-definiteness of $\Lambda^E_{123}$
for inversive distance circle packing metrics
is proved by Guo \cite{Guo} for nonnegative inversive distance and by
the author \cite{X2} for inversive distance in $(-1, +\infty)$ satisfying the structure condition $(\ref{structure condition})$.
The proof we give here simplifies the proof in \cite{Guo,X2}.
\end{remark}

\subsection{Proof of the rigidity for Euclidean inversive distance circle packing}
As the rest of the proof for the rigidity is standard and the same as that in \cite{Guo,L,X2},
we just give a sketch of the proof here.
For more details of the proof, see \cite{Guo,L,X2}.

By Lemma \ref{symmetry Euclidean} and Proposition \ref{simply connect of admi space with weight},
we can define the following function
$$F_{ijk}(u_i, u_j, u_k)=\int_{(\overline{u}_i, \overline{u}_j, \overline{u}_k)}^{(u_i, u_j, u_k)}\theta_idu_i+\theta_jdu_j+\theta_kdu_k$$
on $\ln(\Omega^E_{ijk}(\eta))$,
which is locally concave by Lemma \ref{negativity of Lambda}.
Recall the following definition and extension theorem of Luo \cite{L}.

\begin{definition}
A differential 1-form $w=\sum_{i=1}^n a_i(x)dx^i$ in an open set $U\subset \mathbb{R}^n$ is said to be continuous if
each $a_i(x)$ is continuous on $U$.  A continuous differential 1-form $w$ is called closed if $\int_{\partial \tau}w=0$ for each
triangle $\tau\subset U$.
\end{definition}

\begin{theorem}[\cite{L}, Corollary 2.6]\label{Luo's convex extention}
Suppose $X\subset \mathbb{R}^n$ is an open convex set and $A\subset X$ is an open subset of $X$ bounded by a real analytic
codimension-1 submanifold in $X$. If $w=\sum_{i=1}^na_i(x)dx_i$ is a continuous closed 1-form on $A$ so that
$F(x)=\int_a^x w$ is locally convex on $A$ and each $a_i$ can be extended continuous to $X$ by constant functions to a
function $\widetilde{a}_i$ on $X$, then  $\widetilde{F}(x)=\int_a^x\sum_{i=1}^n\widetilde{a}_i(x)dx_i$ is a $C^1$-smooth
convex function on $X$ extending $F$.
\end{theorem}

By Proposition \ref{simply connect of admi space with weight},
Lemma \ref{continuity of Euclidan inner angle} and Theorem \ref{Luo's convex extention},
$F_{ijk}$ could be extended to be a $C^1$-smooth concave function
$$\widetilde{F}_{ijk}(u_i, u_j, u_k)=\int_{(\overline{u}_i, \overline{u}_j, \overline{u}_k)}^{(u_i, u_j, u_k)}\widetilde{\theta}_idu_i+\widetilde{\theta}_jdu_j+\widetilde{\theta}_kdu_k$$
defined on $\mathbb{R}^3$.
Using $\widetilde{F}_{ijk}$, we can further define the following $C^1$ convex function $\widetilde{F}$ on $\mathbb{R}^{|V|}$
$$\widetilde{F}(u_1, \cdots, u_{|V|})=2\pi\sum_{i\in V}u_i-\sum_{\{ijk\}\in F}\widetilde{F}_{ijk}(u_i, u_j, u_k),$$
which has gradient $\nabla_{u_i} \widetilde{F}=2\pi-\sum_{\{ijk\}\in F}\widetilde{\theta}_{i}^{jk}=\widetilde{K}_i$,
where $\widetilde{K}_i$ is a continuous extension of $K_i$.
Then the global rigidity of $K$ on the admissible space of Euclidean inversive distance circle packing metrics for $(M, \mathcal{T}, I)$
follows from the convexity of $\widetilde{F}$ and the null space of $\Lambda^E_{ijk}(\eta)$ is
$\{t(1,1,1)|t\in \mathbb{R}\}$. \qed

\section{Rigidity of hyperbolic inversive distance circle packing}\label{section 3}
\subsection{The admissible space of hyperbolic inversive distance circle packing metrics for a single triangle}

Similar to the Euclidean case, we can define the admissible space of hyperbolic inversive
distance circle packing metrics for a triangle $\sigma=\{123\}\in F$.
Given a weight $\eta$ on the edge set $E_\sigma$ of $\{123\}$ satisfying the structure condition (\ref{structure condition}),
the admissible space $\Omega^H_{123}(\eta)$ of hyperbolic inversive distance circle packing metrics for the triangle $\sigma=\{123\}$
is defined to be the set of hyperbolic inversive distance circle packing metrics $(r_1, r_2, r_3)\in \mathbb{R}^3_{>0}$ such that
the triangle with edge lengths given by (\ref{hyperbolic length})
exists in 2-dimensional hyperbolic space $\mathbb{H}^2$.

Set
\begin{equation*}
\begin{aligned}
G_-(l)=\left(
         \begin{array}{cccc}
           0 & 1 & 1 & 1 \\
           1 & 0 & -\cosh l_{12} & -\cosh l_{13} \\
           1 & -\cosh l_{12} & 0 & -\cosh l_{23} \\
           1 & -\cosh l_{13} & -\cosh l_{23} & 0 \\
         \end{array}
       \right).
\end{aligned}
\end{equation*}
Recall the following result characterizing
nondegeneracy of a hyperbolic triangle $\{123\}$ with positive edge lengths $l_{12}, l_{13}, l_{23}$.

\begin{lemma}[\cite{TW}, Proposition 2.4.1]
A triangle with positive edge lengths $l_{12}, l_{13}, l_{23}$ exists in $\mathbb{H}^2$ if and only
if $\det G_-(l)<0$.
\end{lemma}

\begin{remark}
By direct calculations, we have
$$\det G_-(l)=-4\sinh s\sinh (s-l_{12})\sinh (s-l_{13})\sinh (s-l_{23}),$$
where $s=\frac{1}{2}(l_{12}+l_{13}+l_{23})$ is the semiperimeter.
This implies that $\det G_-(l)<0$ is equivalent to the triangle inequalities.
This was also observed by Guo in \cite{Guo}.
Similar to the Euclidean case, this approach has the advantage that we just need
one inequality to characterize the nondegeneracy instead of three triangle inequalities.
Furthermore, this approach could be generalized to high dimensions \cite{TW}.
\end{remark}

For simplicity, we set
\begin{equation*}
\begin{aligned}
C_i=\cosh r_i, S_i=\sinh r_i.
\end{aligned}
\end{equation*}
Submitting the definition of $l_{ij}$ (\ref{hyperbolic length}) into $G_-(l)$, we have
\begin{equation*}
\begin{aligned}
-\det G_-(l)
=&2S_1^2S_2^2S_3^2(1+I_1I_2I_3)\\
&+S_1^2S_2^2(1-I_3^2)+S_1^2S_3^2(1-I_2^2)+S_2^2S_3^2(1-I_1^2)\\
&+2C_2C_3S_1^2S_2S_3\gamma_1+2C_1C_3S_1S_2^2S_3\gamma_2+2C_1C_2S_1S_2S_3^2\gamma_3.
\end{aligned}
\end{equation*}
Set
\begin{equation*}
\begin{aligned}
\kappa_i=\coth r_i,
\end{aligned}
 \end{equation*}
then
\begin{equation*}
\begin{aligned}
&-\det G_-(l)
=S_1^2S_2^2S_3^2Q,
\end{aligned}
\end{equation*}
where
\begin{equation*}
\begin{aligned}
Q=&\kappa_1^2(1-I_1^2)+\kappa_2^2(1-I_2^2)+\kappa_3^2(1-I_3^2)+2\kappa_1\kappa_2\gamma_3+2\kappa_1\kappa_3\gamma_2+2\kappa_2\kappa_3\gamma_1\\
  &+I_1^2+I_2^2+I_3^2+2I_1I_2I_3-1.
\end{aligned}
 \end{equation*}
Then we have the following criterion of nondegeneracy for hyperbolic triangles.
\begin{lemma}[\cite{Guo,X2}]\label{nondegenerate lemma hyperbolic}
A hyperbolic triangle $\{123\}$ with edge lengths $l_{12}, l_{13}, l_{23}$
given by $(\ref{hyperbolic length})$ exists in $\mathbb{H}^2$ if and only if $Q>0$.
\end{lemma}

Similar to the Euclidean case, set
\begin{equation*}
\begin{aligned}
h_1=&\kappa_1(1-I_1^2)+\kappa_2\gamma_3+\kappa_3\gamma_2,\\
h_2=&\kappa_2(1-I_2^2)+\kappa_1\gamma_3+\kappa_3\gamma_1,\\
h_3=&\kappa_3(1-I_3^2)+\kappa_1\gamma_2+\kappa_2\gamma_1,
\end{aligned}
 \end{equation*}
we have
\begin{equation*}
\begin{aligned}
Q=&\kappa_1h_1+\kappa_2h_2+\kappa_3h_3 +I_1^2+I_2^2+I_3^2+2I_1I_2I_3-1.
\end{aligned}
 \end{equation*}
Then $(r_1, r_2, r_3)\in \mathbb{R}^3_{>0}$ is a degenerate hyperbolic inversive distance circle packing metric
for a single triangle $\{123\}$ if and only if
\begin{equation}\label{Q leq 0 hyperbolic}
\begin{aligned}
Q=\kappa_1h_1+\kappa_2h_2+\kappa_3h_3 +I_1^2+I_2^2+I_3^2+2I_1I_2I_3-1\leq 0.
\end{aligned}
 \end{equation}
If $I_1, I_2, I_3\in (-1, 1]$, we have
$h_i\geq 1-I_i^2+\gamma_j+\gamma_k$ by $\kappa_i=\coth r_i>1$, which implies
\begin{equation*}
\begin{aligned}
Q\geq& (1-I_1^2)+\gamma_3+\gamma_2+(1-I_2^2)+\gamma_3+\gamma_1+\\
&(1-I_3^2)+\gamma_2+\gamma_1+I_1^2+I_2^2+I_3^2+2I_1I_2I_3-1\\
=&2(1+I_1)(1+I_2)(1+I_3)>0
\end{aligned}
\end{equation*}
for any $(r_1,r_2,r_3)\in \mathbb{R}^3_{>0}$.
Therefore, if $(r_1,r_2,r_3)\in \mathbb{R}^3_{>0}$ is a degenerate hyperbolic inversive distance circle packing metric for a triangle $\{123\}$,
at least one of $I_1, I_2, I_3$ is strictly larger than $1$,
which implies $I_1^2+I_2^2+I_3^2+2I_1I_2I_3-1>0$ by the proof of Lemma \ref{positive Delta Euclidean}.

Therefore, if $(r_1,r_2,r_3)\in \mathbb{R}^3_{>0}$ is a degenerate inversive distance circle packing metric for a triangle $\{123\}$,
we have
$$\kappa_1h_1+\kappa_2h_2+\kappa_3h_3< 0$$
by (\ref{Q leq 0 hyperbolic}),
which implies one of the following two cases happens.
\begin{description}
  \item[(1)] At least one of $h_1,h_2,h_3$ is zero;
  \item[(2)] None of $h_1,h_2,h_3$ is zero.
\end{description}
Similar to the Euclidean case, we can prove that case (1) never happens. Furthermore,
we can prove that only one of $h_i,h_j,h_k$ is negative and the others are positive in case (2).

Similar to the Euclidean case, we can rewrite $Q\leq 0$ as a quadratic inequality of $\kappa_i$
\begin{equation*}
A_i\kappa_i^2+B_i\kappa_i+C_i\geq 0,
\end{equation*}
where
\begin{equation}\label{A1,B1,C1 hyperbolic}
\begin{aligned}
A_i=&I_i^2-1,\\
B_i=&-2(\kappa_j\gamma_k+\kappa_k\gamma_j)\leq 0, \ \ \{i,j,k\}=\{1,2,3\},\\
C_i=&\kappa_j^2(I_j^2-1)+\kappa_k(I_k^2-1)-2\kappa_j\kappa_k\gamma_i\\
&-(I_1^2+I_2^2+I_3^2+2I_1I_2I_3-1)
\end{aligned}
 \end{equation}
with $\{i,j,k\}=\{1,2,3\}$.
By direct calculations, we have the determinant $\Delta_i=B_i^2-4A_iC_i$ is given by
\begin{equation}\label{Delta-i hyperbolic}
\begin{aligned}
\Delta_i
=&4(\kappa_j^2+\kappa_k^2+2\kappa_j\kappa_kI_i)(I_1^2+I_2^2+I_3^2+2I_1I_2I_3-1)\\
 &+4(I_i^2-1)(I_1^2+I_2^2+I_3^2+2I_1I_2I_3-1).
\end{aligned}
 \end{equation}
Similar to the Euclidean case, we have the following results.

\begin{lemma}\label{positive Delta hyperbolic}
If $I_i>1$ and the structure condition $(\ref{structure condition})$ is satisfied, then $\Delta_i>0$.
\end{lemma}

\begin{lemma}\label{lemma none of h1,h2,h3 is 0 hyperbolic}
Suppose $(r_1, r_2, r_3)\in \mathbb{R}^3_{>0}$ is a degenerate hyperbolic inversive distance circle packing
metric for a triangle $\{123\}$
with a weight $\eta: E_\sigma\rightarrow (-1, +\infty)$
satisfying the structure condition $(\ref{structure condition})$, then none of $h_1, h_2, h_3$ is zero.
\end{lemma}

\begin{lemma}\label{lemma no h1<0 h2<0 hyperbolic}
Suppose $\{123\}\in F$ is a triangle
with a weight $\eta: E_\sigma\rightarrow (-1, +\infty)$
satisfying the structure condition $(\ref{structure condition})$
and $(r_1, r_2, r_3)\in \mathbb{R}^3_{>0}$. Then there exists no subset $\{i,j\}\subset \{1,2,3\}$ such that
$h_i<0$ and $h_j<0$.
\end{lemma}

\begin{proposition}[\cite{Guo,X2}]\label{simply connect of admi space with weight hyperbolic}
Suppose $\sigma=\{123\}\in F$ is a triangle in $(M, \mathcal{T})$ with a weight $\eta: E_\sigma\rightarrow (-1, +\infty)$
satisfying the structure condition $(\ref{structure condition})$.  Then the admissible
space $\Omega^H_{123}(\eta)$ of hyperbolic inversive distance circle packing metrics $(r_1, r_2, r_3)\in \mathbb{R}^3_{>0}$
is nonempty and simply connected.
Furthermore, the set of degenerate inversive distance circle packing metric is a disjoint union
$\cup_{i\in P}V_i,$
where $P=\{i\in \{1,2,3\}|I_i>1\}$ and
$$V_i=\left\{(r_1, r_2, r_3)\in \mathbb{R}^3_{>0}|\kappa_i\geq \frac{-B_i+\sqrt{\Delta_i}}{2A_i}\right\}$$
is bounded by an analytic graph on $\mathbb{R}^2_{>0}$
with $A_i, B_i, C_i, \Delta_i$ given by $(\ref{A1,B1,C1 hyperbolic}) (\ref{Delta-i hyperbolic})$.
\end{proposition}

Lemma \ref{positive Delta hyperbolic}, Lemma \ref{lemma none of h1,h2,h3 is 0 hyperbolic},
Lemma \ref{lemma no h1<0 h2<0 hyperbolic} and Proposition \ref{simply connect of admi space with weight hyperbolic}
could be proved similarly to that of Lemma \ref{positive Delta Euclidean}, Lemma \ref{lemma none of h1,h2,h3 is 0},
Lemma \ref{lemma no h1<0 h2<0 Euclidean} and Proposition \ref{simply connect of admi space with weight}
by repeating the proof line by line.
We omit the details of the proof here. Similar to Remark \ref{Vi hi<0, hj,hk>0 Euclidean}, we have the following remark.

\begin{remark}\label{Vi hi<0, hj,hk>0 hyperbolic}
If $(r_1,r_2,r_3)\in V_i$ is a degenerate hyperbolic inversive distance circle packing metric
for a triangle $\{123\}$
with a weight $\eta: E_\sigma\rightarrow (-1, +\infty)$
satisfying the structure condition $(\ref{structure condition})$, then $h_i<0$, $h_j>0$, $h_k>0$, where $\{i,j,k\}=\{1,2,3\}$.
\end{remark}

Similar to the Euclidean case, the inner angles of a hyperbolic triangle could be extended by constants to be globally
defined continuous functions.
\begin{lemma}[\cite{L,X2}]\label{continuity of hyperbolic inner angle}
Suppose $\{123\}\in F$ is a triangle
with a weight $\eta: E_\sigma\rightarrow (-1, +\infty)$
satisfying the structure condition $(\ref{structure condition})$.
Then the functions $\theta_1, \theta_2, \theta_3$ defined for $(r_1, r_2, r_3)\in \Omega^H_{123}(\eta)$
could be extended by constants to be continuous functions
$\widetilde{\theta}_1, \widetilde{\theta}_2, \widetilde{\theta}_3$ defined on $\mathbb{R}^3_{>0}$.
\end{lemma}
\proof
Suppose $(r_1,r_2,r_3)\in \Omega^H_{123}(\eta)$ tends to a point $(\overline{r}_1,\overline{r}_2,\overline{r}_3)\in \partial V_i$.
By direct calculations, we have
\begin{equation*}
\begin{aligned}
-\det G_-(l)
=&4\sinh s\sinh (s-l_{ij})\sinh (s-l_{ik})\sinh (s-l_{jk})\\
=&(\cosh (l_{jk}+l_{ik})-\cosh l_{ij})(\cosh l_{ij}-\cosh (l_{jk}-l_{ik}))\\
=&(\cosh^2 l_{jk}-1)(\cosh l_{ik}^2-1)-(\cosh l_{jk}\cosh l_{ik}-\cosh l_{ij})^2\\
=&\sinh^2l_{jk}\sinh^2l_{ik} -\sinh^2l_{jk}\sinh^2l_{ik} \cos^2\theta_k\\
=&\sinh^2l_{jk}\sinh^2l_{ik}\sin^2\theta_k,
\end{aligned}
\end{equation*}
where $\{i,j,k\}=\{1,2,3\}$.
As $(r_1,r_2,r_3)\in \Omega^H_{123}(\eta)$ tends to $(\overline{r}_1,\overline{r}_2,\overline{r}_3)\in \partial V_i$,
we have $\det G_-(l)\rightarrow 0$, which implies $\theta_1, \theta_2,\theta_3\rightarrow 0$ or $\pi$.

By Lemma \ref{symmetry hyperbolic}, we have
\begin{equation*}
\begin{aligned}
\frac{\partial \theta_j}{\partial u_i}
=\frac{S_i^2S_j^2S_k}{2\widetilde{A}_{123}\sinh^2 l_{ij}}[\kappa_k(1-I_k^2)+\kappa_i\gamma_j+\kappa_j\gamma_i]
=\frac{S_i^2S_j^2S_kh_k}{2\widetilde{A}_{123}\sinh^2 l_{ij}},
\end{aligned}
\end{equation*}
where $u_i=\ln \tanh\frac{r_i}{2}$ and $\widetilde{A}_{123}=\frac{1}{2}\sinh l_{ik}\sinh l_{ij}\sin \theta_i$.
Note that for $(\overline{r}_1,\overline{r}_2,\overline{r}_3)\in \partial V_i$, we have $h_k>0$ by Remark \ref{Vi hi<0, hj,hk>0 hyperbolic}.
Therefore, for  $(r_1,r_2,r_3)\in \Omega^H_{123}(\eta)$ sufficiently close to $(\overline{r}_1,\overline{r}_2,\overline{r}_3)\in \partial V_i$,
we have $\frac{\partial \theta_j}{\partial u_i}>0$,
which implies $\theta_j,\theta_k\rightarrow 0$ as $(r_1,r_2,r_3)\rightarrow (\overline{r}_1,\overline{r}_2,\overline{r}_3)\in \partial V_i$.

Furthermore, we have the following formula \cite{V} for the area $A_{123}$ of the hyperbolic triangle $\{123\}$
\begin{equation*}
\begin{aligned}
\tanh^2\frac{A_{123}}{4}
=&\tanh \frac{p}{2}\tanh \frac{p-l_{12}}{2}\tanh \frac{p-l_{13}}{2}\tanh \frac{p-l_{23}}{2}\\
=&\frac{-\det G_-(l)}{64\cosh^2 \frac{p}{2}\cosh^2 \frac{p-l_{12}}{2}\cosh^2 \frac{p-l_{13}}{2}\cosh^2 \frac{p-l_{23}}{2}},
\end{aligned}
\end{equation*}
where $p=\frac{1}{2}(l_{12}+l_{13}+l_{23})$.
This implies $A_{123}\rightarrow 0$ as $(r_1,r_2,r_3)\rightarrow (\overline{r}_1,\overline{r}_2,\overline{r}_3)\in \partial V_i$.
Further note that
$A_{123}=\pi-\theta_{1}-\theta_{2}-\theta_{3}$ and $\theta_j, \theta_k\rightarrow 0$,
we have $\theta_{i}\rightarrow \pi$ as $(r_1,r_2,r_3)\rightarrow (\overline{r}_1,\overline{r}_2,\overline{r}_3)\in V_i$.

Therefore, we can extend $\theta_1, \theta_2,\theta_3$ defined on $\Omega^H_{123}(\eta)$ to be continuous functions defined on $\mathbb{R}^3_{>0}$
by setting
\begin{equation*}
\begin{aligned}
\widetilde{\theta}_i(r_1,r_2,r_3)=\left\{
                       \begin{array}{ll}
                         \theta_i, & \hbox{if $(r_1,r_2,r_3)\in \Omega^H_{123}(\eta)$;} \\
                         \pi, & \hbox{if $(r_1,r_2,r_3)\in V_i$;} \\
                         0, & \hbox{otherwise.}
                       \end{array}
                     \right.
\end{aligned}
\end{equation*}
\qed

Similar to the Euclidean case, we can define the following $6$-dimensional parameterized admissible space
\begin{equation*}
\begin{aligned}
\Omega^H_{123}=\cup_{\eta\in \Gamma}\Omega^H_{123}(\eta).
\end{aligned}
\end{equation*}

\begin{lemma}\label{connectivity of Omega123 hyperbolic}
$\Omega^H_{123}$ is connected.
\end{lemma}
The proof of Lemma \ref{connectivity of Omega123 hyperbolic} is the same as that of Lemma \ref{connectivity of Omega123 Euclidean},
we omit the proof here.

\subsection{Negative definiteness of the Jacobian matrix}

Set $u_i=\ln \tanh \frac{r_i}{2}$.
The following result on the matrix $\Lambda^H_{123}=\frac{\partial(\theta_1, \theta_2, \theta_3)}{\partial (u_1, u_2, u_3)}$
is known.

\begin{lemma}[\cite{CL,Guo,T1,X2,Z}]\label{symmetry hyperbolic}
Suppose $(r_1,r_2,r_3)\in \mathbb{R}^3_{>0}$ is a nondegenerate hyperbolic inversive distance circle packing metric
for a triangle $\{123\}$ with a weight $\eta$ satisfying the structure condition $(\ref{structure condition})$, then
\begin{equation}\label{derivative of theta hyperbolic}
\frac{\partial \theta_i}{\partial u_j}=\frac{\partial \theta_j}{\partial u_i}
=\frac{S_i^2S_j^2S_k}{2\widetilde{A}_{123}\sinh^2 l_{ij}}\left[\kappa_k(1-I_k^2)+\kappa_j\gamma_{i}+\kappa_i\gamma_{j}\right],
\end{equation}
where $\widetilde{A}=\frac{1}{2}\sinh l_{ik}\sinh l_{ij}\sin \theta_i$ and $\{i,j,k\}=\{1,2,3\}$.
Specially, for $\eta=(I_1, I_2, I_3)=(1,1,1)$, the matrix
$\Lambda^H_{123}=\frac{\partial(\theta_1, \theta_2, \theta_3)}{\partial (u_1, u_2, u_3)}$
is negative definite at $(r_1,r_2, r_3)=(1,1,1)$.
\end{lemma}

\begin{remark}
By $(\ref{derivative of theta hyperbolic})$ and Remark \ref{Vi hi<0, hj,hk>0 hyperbolic},
if $\eta\in \Gamma$ and $(r_1,r_2,r_3)\in \Omega^H_{123}(\eta)$
tends to a point $(\overline{r}_1, \overline{r}_2, \overline{r}_3)\in \partial V_i$,
we have $\frac{\partial \theta_i}{\partial u_j}\rightarrow +\infty$, $\frac{\partial \theta_i}{\partial u_k}\rightarrow +\infty$.
Recall the following formula obtained in Proposition 9 of \cite{GT}
\begin{equation*}
\frac{\partial A_{123}}{\partial u_i}=\frac{\partial \theta_j}{\partial u_i}(\cosh l_{ij}-1)+\frac{\partial \theta_k}{\partial u_i}(\cosh l_{ik}-1),
\end{equation*}
we have $\frac{\partial A_{123}}{\partial u_i}\rightarrow +\infty$, which implies
\begin{equation*}
\frac{\partial \theta_i}{\partial u_i}=-\frac{\partial A_{123}}{\partial u_i}-\frac{\partial \theta_j}{\partial u_i}-\frac{\partial \theta_k}{\partial u_i}
\rightarrow -\infty.
\end{equation*}
\end{remark}

\begin{lemma}\label{negativity of Lambda hyperbolic}
Suppose $r=(r_1, r_2, r_3)\in \mathbb{R}^3_{>0}$ is a nondegenerate hyperbolic inversive distance circle packing metric
for a triangle $\{123\}$ with a weight $\eta$ satisfying the structure condition $(\ref{structure condition})$, then
the matrix $\Lambda^H_{123}(\eta)=\frac{\partial(\theta_1, \theta_2, \theta_3)}{\partial (u_1, u_2, u_3)}$
is negative definite.
\end{lemma}
\proof
By the chain rules, we have
\begin{equation*}
\begin{aligned}
\Lambda^H_{123}
=\frac{\partial(\theta_1, \theta_2, \theta_3)}{\partial (l_{23}, l_{13}, l_{12})}
\cdot
\frac{\partial (l_{23},l_{13},l_{12})}{\partial (u_1, u_2, u_3)}.
\end{aligned}
\end{equation*}
By direct calculations, we have
\begin{equation*}
\begin{aligned}
\det \frac{\partial (l_{23},l_{13},l_{12})}{\partial (u_1, u_2, u_3)}
=&\frac{S_1S_2S_3}{\sinh l_{12}\sinh l_{13}\sinh l_{23}}\cdot\\
&\big[2C_1C_2C_3S_1S_2S_3(1+I_1I_2I_3)+C_1S_1\gamma_1(C_2^2S_3^2 +C_3^2S_2^2)\\
 &+C_2S_2\gamma_2(C_1^2S_3^2 +C_3^2S_1^2)+C_3S_3\gamma_3(C_1^2S_2^2+C_2^2S_1^2)\big]\\
\geq &\frac{2C_1C_2C_3S_1^2S_2^2S_3^2}{\sinh l_{12}\sinh l_{13}\sinh l_{23}}(1+I_1I_2I_3+\gamma_1+\gamma_2+\gamma_3)\\
=&\frac{2C_1C_2C_3S_1^2S_2^2S_3^2}{\sinh l_{12}\sinh l_{13}\sinh l_{23}}(1+I_1)(1+I_2)(1+I_3)>0,
\end{aligned}
\end{equation*}
which implies the matrix $\frac{\partial (l_{23},l_{13},l_{12})}{\partial (u_1, u_2, u_3)}$
is nondegenerate. Therefore,  the rank of $\Lambda^H_{123}$ is the same as that of the matrix
$\frac{\partial(\theta_1, \theta_2, \theta_3)}{\partial (l_{23}, l_{13}, l_{12})}$,
which is 3 for nondegenerate hyperbolic triangles.

Taking $\Lambda^H_{123}$ as a matrix-valued function defined on $\Omega^H_{123}$, then
the three nonzero eigenvalues of $\Lambda^H_{123}$ are continuous functions of $(\eta, r)\in \Omega^H_{123}$.
By the connectivity of $\Omega^H_{123}$ in Lemma \ref{connectivity of Omega123 hyperbolic},
the three nonzero eigenvalues do not change sign on $\Omega^H_{123}$.
Note that $\Lambda^H_{123}(\eta_0)$ is negative definite at
$(r_1, r_2, r_3)=(1,1,1)$ for $\eta_0=(1,1,1)$ by Lemma \ref{negativity of Lambda hyperbolic},
we have
$\Lambda^H_{123}(\eta)$ is negative definite for any $\eta\in \Gamma$.\qed

\subsection{Proof of the rigidity for hyperbolic inversive distance circle packing metrics}
Similar to the Euclidean case, we just sketch the proof
of rigidity for hyperbolic inversive distance circle packing here.
For more details of the proof, see \cite{Guo,L,X2}.

By Lemma \ref{symmetry hyperbolic} and Proposition \ref{simply connect of admi space with weight hyperbolic},
we can define the following function
$$F_{ijk}(u_i, u_j, u_k)=\int_{(\overline{u}_i, \overline{u}_j, \overline{u}_k)}^{(u_i, u_j, u_k)}\theta_idu_i+\theta_jdu_j+\theta_kdu_k$$
on the image of $\Omega^H_{ijk}(\eta)$ under the map $u_i=\ln \tanh \frac{r_i}{2}$,
which is locally concave by Lemma \ref{negativity of Lambda hyperbolic}.
By Proposition \ref{simply connect of admi space with weight hyperbolic},
Lemma \ref{continuity of hyperbolic inner angle} and Theorem \ref{Luo's convex extention},
$F_{ijk}$ could be extended to be a $C^1$-smooth concave function
$$\widetilde{F}_{ijk}(u_i, u_j, u_k)=\int_{(\overline{u}_i, \overline{u}_j, \overline{u}_k)}^{(u_i, u_j, u_k)}\widetilde{\theta}_idu_i+\widetilde{\theta}_jdu_j+\widetilde{\theta}_kdu_k$$
defined on $\mathbb{R}^3_{<0}$.
Using $\widetilde{F}_{ijk}$, we can further define the following $C^1$ convex function $\widetilde{F}$ on $\mathbb{R}^{|V|}_{<0}$
$$\widetilde{F}(u_1,\cdots, u_{|V|})=2\pi\sum_{i\in V}u_i-\sum_{\{ijk\}\in F}\widetilde{F}_{ijk}(u_i, u_j, u_k),$$
which has gradient $\nabla_{u_i} \widetilde{F}=2\pi-\sum_{\{ijk\}\in F}\widetilde{\theta}_{i}^{jk}=\widetilde{K}_i$,
where $\widetilde{K}_i$ is a continuous extension of $K_i$.
Then the global injectivity of $K$ on the admissible space of hyperbolic inversive distance circle packing metrics for $(M, \mathcal{T}, I)$
follows from the convexity of $\widetilde{F}$. This is equivalent to the global rigidity of the curvature map $K$. \qed
\\

\textbf{Acknowledgements}\\[8pt]
This work was partially done when the author was visiting the Department of Mathematics, Rutgers University.
The author thanks Professor Feng Luo for invitation and communications.
The research of the author is supported by Hubei Provincial Natural Science Foundation of China under grant no. 2017CFB681,
Fundamental Research Funds for the Central Universities under grant no. 2042018kf0246 and
National Natural Science Foundation of China under grant no. 61772379 and no. 11301402.

(Xu Xu) School of Mathematics and Statistics, Wuhan University, Wuhan 430072, P.R. China

E-mail: xuxu2@whu.edu.cn\\[2pt]


\begin{thebibliography}{50}
\setlength{\itemsep}{-1pt} \small

\bibitem{An1} E. M. Andreev, \emph{Convex polyhedra in Lobachevsky spaces}. (Russian) Mat. Sb. (N.S.) 81 (123) 1970 445-478.

\bibitem{An2}E. M. Andreev, \emph{Convex polyhedra of finite volume in Lobachevsky space}. (Russian) Mat. Sb. (N.S.) 83 (125) 1970 256-260.

\bibitem{BPS} A. Bobenko, U. Pinkall, B. Springborn, \emph{Discrete conformal maps and ideal hyperbolic polyhedra}.  Geom. Topol. 19 (2015), no. 4, 2155-2215.

\bibitem{BB} J.C. Bowers, P.L. Bowers, \emph{Ma-Schlenker c-Octahedra in the $2$-Sphere}.  Discrete Comput. Geom. 60 (2018), no. 1, 9-26.

\bibitem{BBP} J.C. Bowers, P.L. Bowers, K. Pratt, \emph{Rigidity of circle polyhedra in the 2-sphere and of hyperideal polyhedra in hyperbolic $3$-space}, To appear in Trans. Amer. Math. Soc. DOI: https://doi.org/10.1090/tran/7483

\bibitem{BH} P. L. Bowers, M. K. Hurdal, \emph{Planar conformal mappings of piecewise flat surfaces}. Visualization
and mathematics III, 3-34, Math. Vis., Springer, Berlin, 2003.

\bibitem{BS} P. L. Bowers, K. Stephenson, \emph{Uniformizing dessins and Bely\u{i} maps via circle packing}. Mem.
Amer. Math. Soc. 170 (2004), no. 805

\bibitem{CL} B. Chow, F. Luo, \emph{Combinatorial Ricci flows on surfaces}, J. Differential Geometry, 63 (2003), 97-129.


\bibitem{DV} Y. C. de Verdi\`{e}re, \emph{Un principe variationnel pour les empilements de cercles}, Invent. Math. 104(3) (1991) 655-669.

\bibitem{GH} H. Ge, B. Hua, \emph{$3$-dimensional combinatorial Yamabe flow in hyperbolic background geometry}, \href{https://arxiv.org/abs/1805.10643} {arXiv:1805.10643 [math.DG].}

\bibitem{GJ1} H. Ge, W. Jiang, \emph{On the deformation of inversive distance circle packings, I}. To appear in Trans. Amer. Math. Soc, \href{http://arxiv.org/abs/1604.08317}{arXiv:1604.08317 [math.GT]}.

\bibitem{GJ2} H. Ge, W. Jiang, \emph{On the deformation of inversive distance circle packings, II}. J. Funct. Anal.  272  (2017),  no. 9, 3573-3595.

\bibitem{GJ3} H. Ge, W. Jiang, \emph{On the deformation of inversive distance circle packings, III}. J. Funct. Anal.  272  (2017),  no. 9, 3596-3609.

\bibitem{GJS} H. Ge, W. Jiang, L. Shen, \emph{On the deformation of ball packings}, \href{https://arxiv.org/abs/1805.10573} {arXiv:1805.10573 [math.DG].}

%
%
%

\bibitem{GX1}  H. Ge, X. Xu, \emph{On a combinatorial curvature for surfaces with inversive distance circle packing metrics}. J. Funct. Anal. 275 (2018), no. 3, 523-558.

%

\bibitem{G1} D. Glickenstein, \emph{Discrete conformal variations and scalar curvature on piecewise flat two and three dimensional manifolds}, J.Differential Geometry,  87 (2011), 201-238.

\bibitem{G2} D. Glickenstein, \emph{Geometric triangulations and discrete Laplacians on manifolds}, \href{https://arxiv.org/abs/math/0508188} {arXiv:math/0508188 [math.MG].}

\bibitem{GT} D. Glickenstein, J.Thomas,  \emph{Duality structures and discrete conformal variations of piecewise constant curvature surfaces}. Adv. Math. 320 (2017), 250-278.

%

\bibitem{Guo} R. Guo, \emph{Local rigidity of inversive distance circle packing}, Trans. Amer. Math. Soc. 363 (2011) 4757-4776.


\bibitem{HX1} X. He, X. Xu,  \emph{Thurston's sphere packing on $3$-dimensional manifolds, I}, Preprint.

\bibitem{HX2} X. He, X. Xu,  \emph{Thurston's sphere packing on $3$-dimensional manifolds, II}, In preparation.

\bibitem{H} Z.-X. He, \emph{Rigidity of infinite disk patterns}. Ann. of Math. (2) 149 (1999), no. 1, 1-33.

\bibitem{HS} M. K. Hurdal, K. Stephenson, \emph{Discrete conformal methods for cortical brain flattening}. Neuroimage, 45 (2009) S86-S98.

\bibitem{K1} P. Koebe, \emph{Kontaktprobleme der konformen Abbildung}. Ber. S\"{a}chs. Akad. Wiss. Leipzig, Math.- Phys. Kl. 88 (1936), 141-164.


\bibitem{L} F. Luo, \emph{Rigidity of polyhedral surfaces, III}, Geom. Topol. 15 (2011), 2299-2319.


\bibitem{MS} J. Ma, J. Schlenker, \emph{Non-rigidity of spherical inversive distance circle packings}.
            Discrete Comput. Geom. 47 (2012), no. 3, 610-617.

\bibitem{MR} A. Marden, B. Rodin, \emph{On Thurston's formulaton and proof of Andreeev's theorem}. Computational methods and function theory (Valpara\'{\i}so, 1989), 103-116, Lecture Notes in Math., 1435, Springer, Berlin, 1990.

\bibitem{S} K. Stephenson, \emph{Introduction to circle packing: The theory of discrete analytic functions}, Cambridge Univ. Press (2005)

\bibitem{TW}  Y.U. Taylor, C.T. Woodward,  \emph{ $6j$  symbols for $U_q(\mathfrak{sl}_2)$ and non-Euclidean tetrahedra}. Selecta Math. (N.S.)  11  (2005),  no. 3-4, 539-571.

\bibitem{T} J. Thomas, \emph{Conformal variations of piecewise constant two and three dimensional manifolds}, Thesis (Ph.D.), The University of Arizona, 2015, 120 pp.

\bibitem{T1} W. Thurston, \emph{Geometry and topology of $3$-manifolds}, Princeton lecture notes 1976, \href{http://www.msri.org/publications/books/gt3m}{http://www.msri.org/publications/books/gt3m}.

\bibitem{V}E.B. Vinberg, \emph{Geometry. II}, Encyclopaedia of Mathematical Sciences, 29, Springer-Verlag, New York, 1988.

\bibitem{X1} X. Xu, \emph{On the global rigidity of sphere packings on $3$-dimensional manifolds}, To appear in J.Differential Geometry, \href{https://arxiv.org/abs/1611.08835v4} {arXiv:1611.08835v4 [math.GT].}

\bibitem{X2} X. Xu, \emph{Rigidity of inversive distance circle packings revisited}, Adv. Math. 332 (2018), 476-509.

\bibitem{ZG} W. Zeng, X, Gu, \emph{Ricci Flow for Shape Analysis and Surface Registration}, Springer Briefs in Mathematics. Springer, New York (2013).

\bibitem{ZZGLG} M. Zhang, W. Zeng, R. Guo, F. Luo,  G. Gu, \emph{Survey on discrete surface Ricci flow}. J. Comput. Sci. Tech. 30 (2015), no. 3, 598-613.

\bibitem{Z} Z. Zhou, \emph{Circle patterns with obtuse exterior intersection angles}, \href{https://arxiv.org/abs/1703.01768}{arXiv:1703.01768v2 [math.GT].}

\end{thebibliography}
\end{document}